\documentclass{article}

\usepackage{amsmath}
\usepackage{amssymb}
\usepackage{bm}
\usepackage{array}
\usepackage{authblk}

\newtheorem{remark}{Remark}

\newcommand{\imag}{{\rm i}}
\newcommand{\diff}{{\rm d}}
\newcommand{\partialn}{ \rm{\partial }}

\begin{document}


\title{Coupling of paraxial and white-noise approximations of the Helmholtz equation in randomly layered media
}



\author{Austin McDaniel\thanks{Austin.McDaniel@asu.edu}~}
\author{Alex Mahalov\thanks{mahalov@asu.edu}}
\affil{School of Mathematical and Statistical Sciences, Arizona State University, Tempe, AZ 85287 USA}
\date{}
\maketitle 

\begin{abstract}
We study the simultaneous paraxial and white-noise limit of the Helmholtz equation in randomly layered media where the refractive index fluctuations are in the direction of propagation.  We consider the regime in which the wavelength is of the same order as the correlation length of the random fluctuations of the refractive index.  We show that this simultaneous limit can be taken in this regime by introducing into the equation an arbitrarily small regularization parameter.  The corresponding paraxial white-noise approximation that we derive is different from that of the previously studied high-frequency regime.  Since the correlation length of the refractive index fluctuations due to atmospheric turbulence varies substantially, our results are relevant for numerous different propagation scenarios including microwave and radiowave propagation through various regions of the atmosphere.
\end{abstract}


\section{Introduction}

Many technologies that rely on electromagnetic wave propagation through various regions of the 
atmosphere are significantly affected by atmospheric turbulence.  In the troposphere, tropopause, and lower stratosphere, random 
fluctuations in temperature, pressure, and humidity give rise to random refractive index fluctuations, referred to as refractive turbulence (see, e.g., \cite{Frehlich,Strohbehn,Tatarskii, McDaniel2017, McDaniel2017OPEX}).  The ionosphere is another turbulent region, where electron density fluctuations lead to fluctuations in the refractive index.  Electromagnetic waves propagating through the atmosphere are affected by these random refractive index fluctuations, which result in phenomena such as beam spreading, loss of spatial coherence, and intensity fluctuations which are referred to as scintillation.  These effects of turbulence on electromagnetic wave propagation have far-reaching detrimental consequences on applications such as astronomical imaging, free-space optical communication, remote sensing, radar, the Global Positioning System (GPS), and imaging systems such as synthetic aperture radar (SAR).  For example, turbulence limits the angular resolution 
of ground-based telescopes, causes speckles in images \cite{Osborn}, and can increase the link error probability of optical communication systems \cite{Andrews,Zhu}.  Ionospheric 
turbulence affects GPS microwave signals \cite{DeLaB} and causes distortions in 
spaceborne SAR images \cite{Gilman}.  

Any theoretical description of the refractive index fluctuations due to atmospheric turbulence must be statistical in nature since these fluctuations cannot be accounted for by deterministic models.
The spatial structure is commonly described in terms of the spectrum of the fluctuations, which is related 
to the covariance function of the refractive index through the Wiener-Khinchin theorem when the fluctuations are statistically homogeneous.
The classical Kolmogorov theory of turbulence assumes that turbulence is both statistically homogeneous and 
isotropic \cite{Frisch}.  
This theory states that the spectrum of the fluctuations contains an inertial 
subrange corresponding to a range of length scales whose upper and lower bounds are referred to as the 
outer scale and inner scale, respectively.  Kolmogorov deduced that the 
three-dimensional spectrum $\Phi (\kappa)$ of the refractive index fluctuations is proportional to $\kappa ^{-11/3}$ over the inertial subrange (this corresponds to a one-dimensional spectrum with a $-5/3$ power law), a power law referred to as the Kolmogorov spectrum.  Non-Kolmogorov spectra, i.e. those with a power law different from $-11/3$, have also been considered for modeling atmospheric turbulence, see, e.g., \cite{Lukin,Mahalov1998,Toselli,Yao}.  The strength of the random fluctuations in the refractive index due to turbulence is usually measured in terms of the structure constant $C_n ^2$ which is a dimensional parameter that is proportional to the variance of the refractive index.
There has been much work on measuring the values of $C_n ^2$ as well as the inner scale and outer scale in different scenarios and different regions of the atmosphere, see, e.g., \cite{Andrews,Tatarskii,Vasseur}.  The outer scale is often identified with the correlation length of the refractive index fluctuations due to turbulence, which is a fundamental quantity in defining the problem that we study in this article.  Within the atmospheric boundary layer the outer scale is commonly taken as $0.4 h$ where $h$ is height 
above the ground.  In the free atmosphere the correlation length of the refractive index can be highly anisotropic, where typical values may be about $100$ m in the vertical direction 
and a few kilometers in the horizontal directions \cite{Tatarskii}.  In the ionosphere the outer scale of turbulence ranges between $1$ and $10$ kilometers \cite{Gilman}.
 
In free space, each component of time-harmonic electric and magnetic fields satisfies
the scalar Helmholtz equation.  When the medium is 
inhomogeneous as in the atmosphere, the components of these fields do not satisfy this equation exactly (the exact description is, of course, given by Maxwell's equations).  Nevertheless, this equation is often a good model that is commonly used for describing wave propagation in situations where the medium is 
inhomogeneous.  In this article we are interested in the case in which waves propagate mainly in a single direction, which we denote by $r_3$, and the effects of refractive index fluctuations in the direction of propagation.  Letting 
$\bm{r} = (r_1 , r_2 , r_3) \in \mathbb{R} ^3$ denote the spatial vector, we consider the Helmholtz equation for a scalar field $\Psi$:
\begin{equation}
\label{stochastic Helmholtz}
 \mathrm{\Delta} \Psi + k_0 ^2 n^2(r_3) \Psi = 0 
\end{equation}
where $k_0$ is the free space wave number and the refractive index $n$ depends only on the spatial coordinate $r_3$.  The wave number $k_0$ is related to the carrier wavelength $\lambda _0$ by $k_0 = 2 \pi / \lambda _0$.  We take the refractive index to be a random function of $r_3$ (i.e., a stochastic process) to model random fluctuations in the direction of propagation, for example those due to atmospheric turbulence.  In such a case the
equation (\ref{stochastic Helmholtz}) is referred to as the Helmholtz equation in randomly layered media, meaning that the refractive index is a random 
function of a single spatial variable.

In many scenarios of interest waves propagate along a privileged axis that is determined by the input beam or source.  In such a situation the paraxial, or parabolic, approximation of the Helmholtz equation is often made.  This approximation amounts to 
reducing the Helmholtz equation to a Schr\"odinger equation in which the evolution variable represents the spatial variable along the axis of 
propagation.  The significant advantage of this approximation is therefore that it replaces the boundary value problem associated with
the Helmholtz equation with an evolution problem.  The resulting Schr\"odinger equation is referred to as the paraxial (or parabolic) wave equation and is used often to describe propagation in 
both deterministic and random media.  For example, it is often the starting point for investigations of nonlinear optics phenomena \cite{Fibich} and is used to study laser beam
propagation through random media \cite{Andrews}.  Numerical algorithms for the paraxial wave equation have been used extensively to solve long-range radiowave propagation problems \cite{Levy}.

In situations where waves propagate through a random medium like the turbulent atmosphere, along with the paraxial approximation another simplification called the white-noise (or Markov) approximation is also commonly made.  This involves approximating the random 
medium fluctuations in the main direction of propagation by white noise, i.e., delta-correlated fluctuations.  The advantage of making this further approximation is that the resulting equation can be analyzed by using the It\^o calculus.  This
enables the closure of the equations for the moments \cite{Fouque1998} and allows explicit expressions for the fourth-order moments to be 
obtained in the scintillation regime \cite{Garnier2015moment}.  The combined use of the paraxial and white-noise approximations has, for example, been employed to study the stability of the wave energy density \cite{Bal2004}, characterize the random shift of the transverse profile of a beam \cite{Garnier2014}, and analyze the effects of fluctuations in the boundaries of waveguides, a situation which is important for applications in underwater acoustics \cite{Borcea}.

The paraxial and white-noise approximations are justified mathematically by taking the limits that correspond to the respective asymptotic scaling regimes described by these approximations.
Because of the advantages of both of these approximations, the problem of simultaneously taking the respective limits that yield the paraxial and white-noise approximations has been
studied for the Helmholtz equation \cite{Bailly} as well as more general wave propagation models \cite{Garnier2015,Gomez}.
In \cite{Bailly} the model corresponding to the simultaneous paraxial white-noise limit was derived in the case of randomly layered media for the 
high-frequency regime where the wavelength is small compared to the correlation length of the fluctuations in the medium.  The model for this high-frequency regime has also been derived in the case of a three-dimensional random medium \cite{Garnier2009,Garnier2009WM}.  In this article we consider the simultaneous paraxial white-noise limit in a different regime, namely the one 
in which the wavelength is of the same order as the correlation length of the random refractive index fluctuations.  We 
show that in this regime there is a coupling between the paraxial and white-noise limits that results in a model that 
is different from the one derived in the high-frequency regime.  As we explain in the Discussion, our new model can be viewed as a generalization of the previously derived model for the high-frequency regime that extends the range of validity to include the regime in which the wavelength is of the same order as the correlation length of the medium fluctuations.  

We let the refractive index $n$ in (\ref{stochastic Helmholtz}) be modeled by
\begin{equation}
\label{refractive index}
n ^2 (r_3) = 1 + \sigma \nu \bigg( \frac{r_3}{\ell _c} \bigg) 
\end{equation}
where $\nu$ is a mean-zero, stationary stochastic process, $\ell _c$ is the correlation length of the random fluctuations, and $\sigma$ is a dimensionless parameter that determines the strength of the random fluctuations. 
The paraxial approximation of the Helmholtz equation is obtained by first expressing the wave field $\Psi$ as
\begin{equation}
\label{udefinition}
 \Psi (\bm{r}) = {\rm e}^{{\rm i} k_0 r_3} u (\bm{r}) \; .
\end{equation}
Since the first factor is a plane wave with wavenumber $k_0$ that propagates in the $r_3$ direction, if the wave described by the solution $\Psi$ of (\ref{stochastic Helmholtz}) propagates mainly in the $r_3$ direction then the amplitude (or envelope) $u$ will be slowly-varying in $r_3$.  The paraxial approximation amounts to taking advantage of this fact in order to obtain a simpler equation for the amplitude $u$.  Substituting (\ref{udefinition}) into (\ref{stochastic Helmholtz}) gives
the following equation for $u$:
\begin{equation}
\label{uequation}
 \Bigg( \frac{\partial ^2}{\partial r_1 ^2} + \frac{\partial ^2}{\partial r_2 ^2}
 \Bigg) u + \frac{\partial ^2 u }{\partial r_3 ^2} + 2 {\rm i} k_0 \frac{\partial u}{\partial r_3}
 + k_0 ^2 \Big( n ^2 (r_3) - 1 \Big) u = 0 \; .
\end{equation}

In order to describe the propagation scenario that we consider in this article it is useful to express equation (\ref{uequation}) in 
dimensionless form.  Let $L$ be the typical propagation distance in the longitudinal direction and let $L_x$ be a reference length scale in the transverse directions.  The length scale $L_x$ is typically taken to be the width of the propagating beam.  We define the dimensionless spatial coordinates
$x_1, x_2,$ and $z$ by
\begin{align*}
  (x_1, x_2) = \left( \frac{r_1}{L_x} , \frac{r_2}{L_x} \right) \hspace{5pt} , \hspace{5pt} z = \frac{r_3}{L} \; .
\end{align*}
We are interested in the propagation regime in which the wavelength is of the same order as the correlation length of the refractive index fluctuations.  We let $\ell$ be a length scale that characterizes the size of the wavelength and correlation length.
We then define in terms of $\ell$ the dimensionless wave number $k$ and the dimensionless correlation length $l_c$:
\begin{equation*}
k = k_0 \ell \hspace{5pt} , \hspace{5pt} l_c = \frac{ \ell_c }{ \ell } \; .
\end{equation*}
In addition, we let 
$\mathrm{\Delta} _{\bm{x}}$ denote the transverse Laplacian
\begin{equation*}
  \mathrm{\Delta} _{\bm{x}} = \frac{\partialn ^2}{\partialn x_1 ^2} + \frac{\partialn ^2}{\partialn x_2 ^2} \; .
\end{equation*}
Then in terms of these dimensionless quantities equation (\ref{uequation}) becomes
\begin{align}
\label{mathematical model a}
  & {\rm i} \frac{\partialn u }{\partialn z} + \frac{ \ell }{2 L k } \frac{\partialn ^2 u}{\partialn z ^2}
 + \frac{L \ell}{2 L_x ^2 k } \mathrm{\Delta} _{\bm{x}} u + \frac{ L k \sigma }{2 \ell } \nu \left( \frac{L z}{ l_c \ell } \right) u  = 0 \; . 
\end{align}

We now describe the paraxial white-noise regime that we consider.  We assume that the propagation distance $L$ is much larger than the wavelength and the correlation 
length of the medium fluctuations.
We therefore define the dimensionless parameter $\epsilon$ by
\begin{equation*}
  \epsilon ^2 = \ell / L 
\end{equation*}
and assume $\epsilon \ll 1$.  The paraxial white-noise regime is characterized by two further assumptions.  
The first of these concerns the strength of the refractive index fluctuations in terms of the small parameter $\epsilon$.  We assume
$\sigma \sim O( \epsilon )$ and let
\begin{equation*}
  \sigma = \beta \epsilon 
\end{equation*}
where $\beta$ is a dimensionless parameter.
The second assumption can be 
described in terms of the Fresnel number $N_{ \mathrm{F}}$, which is defined as
\begin{equation}
\label{Fresnel number}
  N_{ \mathrm{F}} = \frac{ L_x ^2 }{ L \lambda _0 } = \frac{ L_x ^2 k }{ 2 \pi L \ell } \; .
\end{equation}
In classical optics, the Fresnel number provides a criterion for the validity of the near field and far field approximations.  Roughly speaking, 
$N_{ \mathrm{F}} \gg 1$ corresponds to the near field and Fresnel diffraction, whereas when $N_{ \mathrm{F}} \ll 1$ the diffracted wave is considered to be in the far field and can often be well-described by Fraunhofer diffraction \cite{Born}.
The Fresnel number can also be expressed in terms of the Fresnel length $L_{ \mathrm{F} } = \sqrt{ L \lambda _0 }$ (see, e.g., \cite{Strohbehn}) and the Rayleigh length $L_{ \mathrm{R} } = \pi L_x ^2 / \lambda _0$ (for a Gaussian beam):
\begin{equation*}
  N_{ \mathrm{F}} = \left( \frac{ L_x }{ L_{ \mathrm{F} } } \right) ^2 = \frac{ L_{ \mathrm{R} } }{ \pi L } \; .
\end{equation*}
The Rayleigh length is the distance along the propagation direction of a beam from the beam waist to the place where the area of the cross-section is doubled by diffraction in free space.
We assume
\begin{equation}
\label{Fresnel number condition}
  N_{ \mathrm{F}} \sim O(1)
\end{equation}
with respect to the small dimensionless parameter $\epsilon$.  The condition (\ref{Fresnel number condition}) implies that the aspect ratio of the wave beam is small, i.e., $L_x / L \ll 1$.  This is called the small angle approximation and is the usual constraint for the validity of the paraxial approximation.  Condition (\ref{Fresnel number condition}) therefore also implies that the transverse length scale $L_x$, which is usually identified with the beam width, is large compared to the wavelength $\lambda _0$, i.e., $L_x / \lambda _0 \gg 1$.

Expressed in terms of the small parameter $\epsilon$, equation (\ref{mathematical model a}) becomes
\begin{align}
\label{mathematical model}
  & {\rm i} \frac{\partialn u }{\partialn z} + \frac{ \epsilon ^2 }{2 k } \frac{\partialn ^2 u}{\partialn z ^2}
 + \frac{L \ell}{2 L_x ^2 k } \mathrm{\Delta} _{\bm{x}} u + \frac{ k \beta }{2 \epsilon } \nu \left( \frac{z}{ \epsilon ^2 l_c } \right) u  = 0 \; . 
\end{align}
The paraxial white-noise regime corresponds to small values of the dimensionless parameter $\epsilon$.  To describe this regime it is therefore useful to obtain a simpler, effective equation by taking the limit $\epsilon \rightarrow 0$.
In order to do this we regularize the problem in the following way.  We introduce into 
equation (\ref{mathematical model}) an 
arbitrarily small parameter $\delta > 0$ as follows:
\begin{align}
\label{regularization}
  \frac{\partialn u}{\partialn z} + \epsilon ^2 \left( \delta - \frac{{\rm i} }{2 k } \right) & \frac{\partialn ^2 u}{\partialn z ^2}
 - {\rm i} \frac{L \ell }{2 L_x ^2 k } \mathrm{\Delta} _{\bm{x}} u - {\rm i} \frac{ k \beta }{2 \epsilon } \nu 
     \left( \frac{z}{ \epsilon ^2 l_c } \right) u  = 0 \; .
\end{align}
The need for introducing this parameter $\delta$ and its purpose will become clear in the next section.  To study the 
simultaneous paraxial white-noise limit in this regime we will first take the 
limit $\epsilon \rightarrow 0$ of equation (\ref{regularization}) and then take the limit $\delta \rightarrow 0$ of the resulting equation.  We will 
see that the two limits $\epsilon \rightarrow 0$ and $\delta \rightarrow 0$ do not commute.

As a model for the fluctuations in the refractive index, we take $\nu$ to be a stationary Ornstein-Uhlenbeck 
process \cite{Gardiner}.
We define $\nu$ as the stationary solution of the stochastic differential equation (SDE)
\begin{equation}
\label{nu SDE}
 {\rm d}\nu _z = - \nu _z {\rm d}z + {\rm d} \tilde{W} _z
\end{equation}
where $\tilde{W}$ is a Wiener process.   Equation~\eqref{nu SDE} has a unique stationary distribution and the stationary solution
of \eqref{nu SDE} is the solution with the initial condition distributed according to this stationary distribution.  Defined this way, $\nu$ is a stationary Gaussian process with mean zero and autocovariance function 
\begin{equation*}
	E \big[ \nu _z \hspace{1pt} \nu _{z'} \big] 
	= \frac{1}{2} {\rm e}^{- |z - z'|} \; .
\end{equation*} 
Then, defining the process $\eta$ by
$\eta (z) = \nu \big( z / ( \epsilon ^2 l_c ) \big) $, we have that $\eta$ satisfies the equation
\begin{equation}
\label{fastOU}
 {\rm d}\eta _z = - \frac{1}{\epsilon ^2 l_c } \eta _z {\rm d}z + \frac{1}{\epsilon \sqrt{l_c} } {\rm d}W _z
\end{equation}
where $W$ is the Wiener process given by $W _z = \epsilon \sqrt{l_c} \tilde{W} _{z / ( \epsilon ^2 l_c )}$.   
The autocovariance function of $\eta$ is 
\begin{equation}
\label{covariance}
	E \big[ \eta _z \hspace{1pt} \eta _{z'} \big] 
	= \frac{1}{2} {\rm e}^{- \frac{|z - z'|}{\epsilon ^2 l_c}} \; .
\end{equation} 
It can be seen from (\ref{covariance}) that the process $\eta$ has correlation length $\epsilon ^2 l_c$. 

By letting $ v = \epsilon \big( \partialn u / \partialn z \big) $, the system that we study, composed of equations (\ref{regularization}) and (\ref{fastOU}), can be written as
\begin{subequations}
\label{firstordersystem}
\begin{align}
\label{firstordersystema}
  {\rm d} u (z, \bm{x}) = & \; \frac{1}{\epsilon} v (z, \bm{x} ) {\rm d}z \\
\label{firstordersystemb}
  {\rm d} v (z, \bm{x} ) = & - \frac{1}{\epsilon ^2} \Bigg( \frac{2 k }{2 \delta k - {\rm i}} \Bigg) v (z, \bm{x} ) {\rm d}z \notag \\
 & + \frac{1}{\epsilon} \frac{ L \ell}{ L_x ^2 } \Bigg( \frac{ {\rm i} }{2 \delta k - {\rm i}} \Bigg)  \mathrm{\Delta} _{\bm{x}} u (z, \bm{x} ) {\rm d}z 
   + \frac{1}{ \epsilon ^2} \Bigg( \frac{ k^2 \beta {\rm i} }{2 \delta k - {\rm i}} \Bigg) u (z, \bm{x} ) \eta (z) {\rm d}z \\
\label{firstordersystemc}
  {\rm d} \eta (z) = & - \frac{1}{\epsilon ^2 l_c } \eta (z) {\rm d}z + \frac{1}{\epsilon \sqrt{l_c}} {\rm d}W (z) \; .
\end{align}
\end{subequations}
We will show in the next section that the limiting equation corresponding to the system (\ref{firstordersystem}) as $\epsilon \rightarrow 0$ is
\begin{flalign}
\label{first limiting equation}
  {\rm d} u (z, \bm{x}) = & \hspace{4pt}  \frac{{\rm i}}{4 \pi N_{ \mathrm{F}}} \mathrm{\Delta} _{\bm{x}} u (z, \bm{x}) {\rm d}z \notag \\
  & - \frac{ k^2 \beta ^2 l_c }
                {8 \bigg( 1 + l_c ^{-1} \Big( \delta - \frac{{\rm i}}{2 k} \Big) \bigg) } u (z, \bm{x}) {\rm d}z 
        + {\rm i} \frac{ k \beta \sqrt{l_c} }{2} u (z, \bm{x}) {\rm d}W (z) \; . \hspace{5pt}
\end{flalign}
Taking the regularization parameter $\delta \rightarrow 0$ in the above equation, we get the simultaneous paraxial white-noise approximation
\begin{align}
 \label{limiting equation}
  {\rm d} u (z, \bm{x}) =  & \hspace{4pt}  \frac{{\rm i}}{4 \pi N_{ \mathrm{F}}} \mathrm{\Delta} _{\bm{x}} u (z, \bm{x}) {\rm d}z \notag \\
  & - \frac{ k^2 \beta ^2 l_c }{ 8 \left( 1 - {\rm i} \frac{1}{2 k l_c } \right) }
    u (z, \bm{x}) {\rm d}z + {\rm i} \frac{ k \beta \sqrt{l_c} }{2} u (z, \bm{x}) {\rm d}W (z)  \; .  \hspace{5pt}
\end{align}

\begin{remark}
The stochastic integral in (\ref{limiting equation}) is understood in the It\^o sense.
The paraxial white-noise approximation (\ref{limiting equation}), which we derive in the next section, is different from the widely used model obtained in the high-frequency regime where the wavelength is small compared to the correlation length of the medium fluctuations \cite{Bailly,Garnier2009,Garnier2009WM}.  As we explain in the Discussion, equation (\ref{limiting equation}) can be viewed as a generalization of the previously derived paraxial white-noise approximation for the high-frequency regime that extends the range of validity to include the regime in which the wavelength is of the same order as the correlation length of the medium fluctuations.  
\end{remark}

\section{Derivation}
\label{Section: Derivation}

In this section we derive the paraxial white-noise approximation (\ref{limiting equation}) starting from the system
(\ref{firstordersystem}).
We begin by taking the Fourier transform of (\ref{firstordersystema}) and (\ref{firstordersystemb}) in the lateral spatial coordinates.  
Letting $\bm{x} = (x_1, x_2)$ and $\bm{\kappa} = (\kappa_1 , \kappa_2)$,
we define
\begin{subequations}
\begin{equation}
\label{Fourier transform}
\hat{u} (z , \bm{\kappa}) = \int _{\mathbb{R} ^2} u (z , \bm{x}) {\rm e} ^{{\rm i} \bm{\kappa} \cdot \bm{x}} {\rm d} x_1 {\rm d} x _2
\end{equation}
and 
\begin{equation}
\hat{v} (z , \bm{\kappa}) = \int _{\mathbb{R} ^2} v (z , \bm{x}) {\rm e} ^{{\rm i} \bm{\kappa} \cdot \bm{x}} {\rm d} x_1 {\rm d} x _2 \; .
\end{equation}
\end{subequations}
By transforming (\ref{firstordersystema}) 
and (\ref{firstordersystemb}) we get the SDE system
\begin{align*}
  {\rm d} \hat{u} _z = & \; \frac{1}{\epsilon} \hat{v} _z {\rm d}z \\
  {\rm d} \hat{v} _z = & - \frac{1}{\epsilon ^2} 
  \Bigg( \frac{2 k }{2 \delta k - {\rm i}} \Bigg) \hat{v} _z {\rm d}z \\
   & - \frac{1}{\epsilon} \frac{ L \ell }{ L_x ^2 } 
   \Bigg( \frac{ {\rm i} }{2 \delta  k - {\rm i}} \Bigg) | \bm{\kappa} | ^2 \hat{u} _z {\rm d}z 
   + \frac{1}{ \epsilon ^2} 
 \Bigg( \frac{ k^2  \beta {\rm i} }{2 \delta k - {\rm i}} \Bigg) \hat{u} _z \eta _z {\rm d}z \hspace{20pt} \\
 {\rm d}\eta _z = & - \frac{1}{\epsilon ^2 l_c } \eta _z {\rm d}z + \frac{1}{\epsilon \sqrt{l_c}} {\rm d}W _z \; .
\end{align*}

We study the limit as $\epsilon \rightarrow 0$ of this system by using a multiscale expansion method which makes use 
of the correspondence between SDEs and certain partial differential equations (PDEs).  Associated with a typical SDE is a backward Kolmogorov equation (and its adjoint, a Fokker-Planck equation)
which is a PDE that describes the transition probability density.  We will work with the backward Kolmogorov equation that corresponds to the above SDE system.  In order to write this equation
it is convenient to view the real and imaginary parts of $\hat{u}$ and $\hat{v}$ as separate state variables and consider the (real-valued) system that these 
variables, together with the state variable $\eta$, satisfy.  Let $\hat{u} ^{\mathrm{R}}$ and $\hat{v} ^{\mathrm{R}}$ denote the real parts and $\hat{u} ^{\mathrm{I}}$ and $\hat{v} ^{\mathrm{I}}$ denote
the imaginary parts of $\hat{u}$ and $\hat{v}$, respectively.  Then $\hat{u} ^{\mathrm{R}} , \hat{u} ^{\mathrm{I}} , \hat{v} ^{\mathrm{R}} , \hat{v} ^{\mathrm{I}} ,$ and $\eta$ satisfy the system  
\begin{align*}
  \diff  \hat{u} ^{\mathrm{R}} _z = & \; \frac{1}{\epsilon}  \hat{v} ^{\mathrm{R}} _z \diff z \\
  \diff  \hat{u} ^{\mathrm{I}} _z = & \; \frac{1}{\epsilon}  \hat{v} ^{\mathrm{I}} _z \diff z \\
  \diff  \hat{v} ^{\mathrm{R}} _z = & - \frac{1}{\epsilon ^2} \Bigg( \frac{4 \delta  k ^2 }{4 \delta ^2 k ^2 + 1} \Bigg) \hat{v} ^{\mathrm{R}} _z \diff z 
 + \frac{1}{\epsilon ^2} \Bigg( \frac{ 2 k }{4 \delta ^2 k ^2 + 1} \Bigg) \hat{v} ^{\mathrm{I}} _z \diff z \\
 & + \frac{1}{\epsilon} \frac{ L \ell }{ L_x ^2 } \Bigg( \frac{ 1 }{4 \delta ^2 k ^2 + 1} \Bigg) | \bm{\kappa} | ^2 \hat{u}^{\mathrm{R}} _z \diff z
 + \frac{1}{\epsilon} \frac{ L \ell }{ L_x ^2 } 
           \Bigg( \frac{ 2 \delta k }{4 \delta ^2 k ^2 + 1} \Bigg) | \bm{\kappa} | ^2 \hat{u}^{\mathrm{I}} _z \diff z \\
 & - \frac{ 1 }{ \epsilon ^2} \Bigg( \frac{ k^2  \beta }{4 \delta ^2 k ^2 + 1} \Bigg) \hat{u} ^{\mathrm{R}} _z \eta _z \diff z
 - \frac{1}{ \epsilon ^2} \Bigg( \frac{2 \delta k ^3 \beta }{4 \delta ^2 k ^2 + 1} \Bigg) \hat{u} ^{\mathrm{I}} _z \eta _z \diff z \\
  \diff  \hat{v} ^{\mathrm{I}} _z = & - \frac{1}{\epsilon ^2} \Bigg( \frac{4 \delta k ^2 }{4 \delta ^2 k ^2 + 1} \Bigg) \hat{v} ^{\mathrm{I}} _z \diff z
  - \frac{1}{\epsilon ^2} \Bigg( \frac{ 2 k }{4 \delta ^2 k ^2 + 1} \Bigg) \hat{v} ^{\mathrm{R}} _z \diff z \\
  & + \frac{1}{\epsilon} \frac{ L \ell }{ L_x ^2 } \Bigg( \frac{ 1 }{4 \delta ^2 k ^2 + 1} \Bigg) | \bm{\kappa} | ^2 \hat{u} ^{\mathrm{I}} _z \diff z 
  - \frac{1}{\epsilon} \frac{ L \ell }{ L_x ^2 } 
             \Bigg( \frac{ 2 \delta k }{4 \delta ^2 k ^2 + 1} \Bigg) | \bm{\kappa} | ^2 \hat{u} ^{\mathrm{R}} _z \diff z \\
 & - \frac{ 1 }{ \epsilon ^2} \Bigg( \frac{ k^2  \beta }{4 \delta ^2 k ^2 + 1} \Bigg) \hat{u} ^{\mathrm{I}} _z \eta _z \diff z
 + \frac{ 1 }{ \epsilon ^2} \Bigg( \frac{2 \delta k ^3 \beta }{4 \delta ^2 k ^2 + 1} \Bigg) \hat{u} ^{\mathrm{R}} _z \eta _z \diff z \\
  \diff \eta _z = & - \frac{1}{\epsilon ^2 l_c } \eta _z \diff z + \frac{1}{\epsilon \sqrt{l_c}} \diff W _z \; .
\end{align*}
The backward Kolmogorov equation corresponding to the above SDE system is
\begin{equation}
\label{bKolmogorov}
 \frac{\partialn \rho}{\partialn z} = \Bigg( \frac{1}{\epsilon ^2} \mathcal{L}_0 + \frac{1}{\epsilon} \mathcal{L}_1 \Bigg) \rho
\end{equation}
where
\begin{align*}
 \mathcal{L}_0 = \; & - \frac{4 \delta k ^2 }{4 \delta ^2 k ^2 + 1} \hat{v} ^{\mathrm{R}} \frac{\partialn }{\partialn \hat{v} ^{\mathrm{R}} }
           + \frac{ 2 k }{4 \delta ^2 k ^2 + 1} \hat{v} ^{\mathrm{I}} \frac{\partialn }{\partialn \hat{v} ^{\mathrm{R}} } 
           - \frac{4 \delta k ^2 }{4 \delta ^2 k ^2 + 1} \hat{v} ^{\mathrm{I}} \frac{\partialn }{\partialn \hat{v} ^{\mathrm{I}} } \\
       &   - \frac{ 2 k }{4 \delta ^2 k ^2 + 1} \hat{v} ^{\mathrm{R}} \frac{\partialn }{\partialn \hat{v} ^{\mathrm{I}} } 
         - \frac{ k^2  \beta }{4 \delta ^2 k ^2 + 1} \hat{u} ^{\mathrm{R}}  \eta \frac{\partialn }{\partialn \hat{v} ^{\mathrm{R}} }
         - \frac{2 \delta  k ^3 \beta }{4 \delta ^2 k ^2 + 1} \hat{u} ^{\mathrm{I}}  \eta \frac{\partialn }{\partialn \hat{v} ^{\mathrm{R}} }  \\
       & - \frac{ k^2  \beta }{4 \delta ^2 k ^2 + 1} \hat{u} ^{\mathrm{I}}  \eta \frac{\partialn }{\partialn \hat{v} ^{\mathrm{I}} }
         + \frac{2 \delta k ^3 \beta }{4 \delta ^2 k ^2 + 1} \hat{u} ^{\mathrm{R}}  \eta \frac{\partialn }{\partialn \hat{v} ^{\mathrm{I}} } 
         - \frac{1}{ l_c } \eta \frac{\partialn }{\partialn \eta}  + \frac{1}{2 l_c } \frac{\partialn ^2 }{\partialn \eta ^2} 
\end{align*}
and
\begin{align*}
 \mathcal{L}_1 = & \; \hat{v} ^{\mathrm{R}} \frac{\partialn }{\partialn \hat{u} ^{\mathrm{R}}}  + \hat{v} ^{\mathrm{I}} \frac{\partialn }{\partialn \hat{u} ^{\mathrm{I}}} 
  + \frac{L \ell}{ L_x ^2 } \Bigg( \frac{ 1 }{4 \delta ^2 k ^2 + 1} \Bigg) 
                                  | \bm{\kappa} | ^2 \hat{u}^{\mathrm{R}} \frac{\partialn }{\partialn \hat{v} ^{\mathrm{R}} } \\
 & + \frac{L \ell}{ L_x ^2 } \Bigg( \frac{ 2 \delta k }{4 \delta ^2 k ^2 + 1} \Bigg) | \bm{\kappa} | ^2 \hat{u}^{\mathrm{I}} \frac{\partialn }{\partialn \hat{v} ^{\mathrm{R}} }
  + \frac{L \ell}{ L_x ^2 } \Bigg( \frac{ 1 }{4 \delta ^2 k ^2 + 1} \Bigg) | \bm{\kappa} | ^2 \hat{u} ^{\mathrm{I}} \frac{\partialn }{\partialn \hat{v} ^{\mathrm{I}} } \\
 &  - \frac{L \ell}{ L_x ^2 } \Bigg( \frac{ 2 \delta  k }{4 \delta ^2 k ^2 + 1} \Bigg) 
                                  | \bm{\kappa} | ^2 \hat{u} ^{\mathrm{R}} \frac{\partialn }{\partialn \hat{v} ^{\mathrm{I}} }  \; .
\end{align*}

We use a homogenization method for PDEs to obtain the limiting backward Kolmogorov equation, from which we can read off
the limiting SDE.
We expand $\rho$ in powers of $\epsilon$, i.e.,
\begin{equation}
\label{rhoexpansion}
 \rho = \rho _0 + \epsilon \rho _1 + \epsilon ^2 \rho _2 + ...
\end{equation}
and derive the backward Kolmogorov equation for the limiting density $\rho_0$.  First, substituting
(\ref{rhoexpansion}) in (\ref{bKolmogorov}) and equating terms of the same order in $\epsilon$ gives the equations
\begin{equation}
\label{order1overepsilonsquared}
 O \left( \frac{1}{\epsilon ^2} \right) : \mathcal{L}_0 \rho _0 = 0 \hspace{31pt}
\end{equation}
\begin{equation}
\label{order1overepsilon}
 O \left( \frac{1}{\epsilon} \right) : \mathcal{L}_0 \rho _1 = - \mathcal{L}_1 \rho _0 \hspace{11pt}
\end{equation}
\begin{equation}
\label{order1}
 \hspace{5pt} O(1) : \frac{\partialn \rho _0}{\partialn z} = \mathcal{L}_0 \rho _2 + \mathcal{L}_1 \rho _1 \; .
\end{equation}
Equation (\ref{order1overepsilonsquared}) implies $\rho _0 = \rho _0 \big( z , \hat{u} ^{\mathrm{R}} , \hat{u} ^{\mathrm{I}} \big) $.  Next, we find $\rho _1$ in terms of
$\rho _0 \big( z , \hat{u} ^{\mathrm{R}} , \hat{u} ^{\mathrm{I}} \big) $ by solving equation (\ref{order1overepsilon}), which, since $\rho _0$ does not depend on $\hat{v} ^{\mathrm{R}}$ or
$\hat{v} ^{\mathrm{I}}$, is equivalent to
\begin{equation}
\label{order1overepsilonrewritten}
 \mathcal{L}_0 \rho _1 = - \hat{v} ^{\mathrm{R}} \frac{\partial \rho _0}{\partial \hat{u} ^{\mathrm{R}}} - \hat{v} ^{\mathrm{I}} \frac{\partial \rho _0 }{\partial \hat{u} ^{\mathrm{I}}} \; .
\end{equation}
The solution of (\ref{order1overepsilonrewritten}) is
\begin{align*}
\rho _1 = & - \frac{1}{2  k } \hat{v} ^{\mathrm{R}} \frac{\partialn \rho _0}{\partialn \hat{u} ^{\mathrm{I}}} + \delta \hat{v} ^{\mathrm{I}} \frac{\partialn \rho _0}{\partialn \hat{u} ^{\mathrm{I}}}
+ \frac{ k  \beta l_c }{2} \hat{u} ^{\mathrm{R}} \eta \frac{\partialn \rho _0}{\partialn \hat{u} ^{\mathrm{I}}} \hspace{130pt} \\
 & \hspace{120pt} + \frac{1}{2 k} \hat{v} ^{\mathrm{I}} \frac{\partialn \rho _0}{\partialn \hat{u} ^{\mathrm{R}}} 
      + \delta \hat{v} ^{\mathrm{R}} \frac{\partialn \rho _0}{\partialn \hat{u} ^{\mathrm{R}}}
  - \frac{ k \beta l_c }{2} \hat{u} ^{\mathrm{I}} \eta \frac{\partialn \rho _0}{\partialn \hat{u} ^{\mathrm{R}}} \; .
\end{align*}
Equation (\ref{order1}) implies that the function $ \partialn \rho _0 / \partialn z - \mathcal{L}_1 \rho _1$ is in the range of
the operator $\mathcal{L}_0$.  Thus, for equation (\ref{order1}) to have a solution, by the Fredholm alternative this function must be orthogonal to the nullspace of $\mathcal{L}_0 ^*$ (where $\mathcal{L}_0 ^*$ denotes the adjoint of $\mathcal{L}_0$).  That is,
the following solvability condition must hold:
\begin{equation}
\label{solvability condition}
 \int _{\mathbb{R} ^3} \rho ^{\infty} \left( \frac{\partialn \rho _0}{\partialn z} - \mathcal{L}_1 \rho _1 \right) \diff \hat{v} ^{\mathrm{R}} \diff \hat{v} ^{\mathrm{I}} \diff \eta = 0
\end{equation}
where $\rho ^{\infty}$ solves the equation
\begin{align}
\label{stationaryFP}
\mathcal{L}_0 ^* \rho ^{\infty} = & \;   \frac{4 \delta k ^2 }{4 \delta ^2 k ^2 + 1} 
                        \frac{\partialn }{\partialn \hat{v} ^{\mathrm{R}} } \Big( \hat{v} ^{\mathrm{R}}  \rho ^{\infty} \Big)
 -  \frac{ 2 k }{4 \delta ^2 k ^2 + 1} \hat{v} ^{\mathrm{I}} \frac{\partialn \rho ^{\infty} }{\partialn \hat{v} ^{\mathrm{R}} } \hspace{65pt} \notag \\
& +  \frac{4 \delta k ^2 }{4 \delta ^2 k ^2 + 1} 
                            \frac{\partialn }{\partialn \hat{v} ^{\mathrm{I}} } \Big( \hat{v} ^{\mathrm{I}} \rho ^{\infty} \Big)
  +  \frac{ 2 k }{4 \delta ^2  k ^2 + 1} \hat{v} ^{\mathrm{R}} \frac{\partialn \rho ^{\infty} }{\partialn \hat{v} ^{\mathrm{I}} } \notag \\
 & + \frac{  k^2  \beta }{4 \delta ^2  k ^2 + 1} 
                            \hat{u} ^{\mathrm{R}}  \eta \frac{\partialn \rho ^{\infty}}{\partialn \hat{v} ^{\mathrm{R}} }
 + \frac{2 \delta k ^3 \beta }{4 \delta ^2  k ^2 + 1} 
                            \hat{u} ^{\mathrm{I}}  \eta \frac{\partialn \rho ^{\infty} }{\partialn \hat{v} ^{\mathrm{R}} }  \notag \\
 & + \frac{ k^2  \beta }{4 \delta ^2  k ^2 + 1} 
                            \hat{u} ^{\mathrm{I}}  \eta \frac{\partialn \rho ^{\infty} }{\partialn \hat{v} ^{\mathrm{I}} }
 -  \frac{2 \delta  k ^3 \beta }{4 \delta ^2 k ^2 + 1} 
                             \hat{u} ^{\mathrm{R}}  \eta \frac{\partialn \rho ^{\infty} }{\partialn \hat{v} ^{\mathrm{I}} } \notag \\
 & \hspace{90pt} + \frac{1}{ l_c } \frac{\partialn }{\partialn \eta} \big( \eta \rho ^{\infty} \big) + \frac{1}{2 l_c } \frac{\partialn ^2 \rho ^{\infty} }{\partialn \eta ^2} = 0 \; .
\end{align}
The backward Kolmogorov equation for the limiting density $\rho _0$ will follow from the solvability condition (\ref{solvability condition}).

We now turn to the question of the existence of a nontrivial (i.e., not identically equal to zero) integrable solution of equation (\ref{stationaryFP}).  Note that (\ref{stationaryFP}) is the stationary Fokker-Planck equation 
corresponding to the SDE system
\begin{align}
\label{SDE L0}
  \diff  \hat{v} ^{\mathrm{R}} _z = & -  \frac{4 \delta  k ^2 }{4 \delta ^2 k ^2 + 1} \hat{v} ^{\mathrm{R}} _z \diff z 
 + \frac{ 2 k }{4 \delta ^2  k ^2 + 1} \hat{v} ^{\mathrm{I}} _z \diff z \notag \hspace{40pt} \\
 & - \frac{ k^2  \beta }{4 \delta ^2  k ^2 + 1}  \hat{u} ^{\mathrm{R}} _z \eta _z \diff z 
  - \frac{2 \delta k ^3 \beta }{4 \delta ^2  k ^2 + 1} \hat{u} ^{\mathrm{I}} _z \eta _z \diff z \notag \\
  \diff  \hat{v} ^{\mathrm{I}} _z = & - \frac{4 \delta k ^2 }{4 \delta ^2  k ^2 + 1} \hat{v} ^{\mathrm{I}} _z \diff z
  - \frac{ 2 k }{4 \delta ^2 k ^2 + 1} \hat{v} ^{\mathrm{R}} _z \diff z \notag \\
 & - \frac{ k^2  \beta }{4 \delta ^2 k ^2 + 1} \hat{u} ^{\mathrm{I}} _z \eta _z \diff z 
  + \frac{2 \delta k ^3 \beta }{4 \delta ^2 k ^2 + 1} \hat{u} ^{\mathrm{R}} _z \eta _z \diff z  \notag \\
  \diff \eta _z = & - \frac{1}{ l_c } \eta _z \diff z + \frac{1}{ \sqrt{ l_c } } \diff W _z \; . 
\end{align}
Thus, the statement that an integrable nontrivial solution of (\ref{stationaryFP}) exists is equivalent to the statement that there exists a stationary 
density corresponding to the SDE system (\ref{SDE L0}).
Letting $\bm{V} = \big( \hat{v} ^{\mathrm{R}} , \hat{v} ^{\mathrm{I}} , \eta \big) ^{\mathrm T} $ (where T denotes transpose), the system (\ref{SDE L0}) can be written as
\begin{equation*}
\diff \bm{V} _z = - {\bm \gamma} \bm{V} _z \diff z + \bm{D} \diff W _z
\end{equation*}
where
\begin{equation*} \bm{ \gamma } =
\left[\rule{0cm}{2.3cm}\right. \arraycolsep=12pt \def\arraystretch{1.5} \begin{array}{ccc}
  \frac{4 \delta k ^2 }{4 \delta ^2 k ^2 + 1}  &
   - \frac{ 2 k }{4 \delta ^2 k ^2 + 1} & 
    \frac{ k^2  \beta }{4 \delta ^2 k ^2 + 1} \hat{u} ^{\mathrm{R}} 
   + \frac{2 \delta k ^3 \beta }{4 \delta ^2 k ^2 + 1} \hat{u} ^{\mathrm{I}} \\[1.3em]
    \frac{ 2 k }{4 \delta ^2 k ^2 + 1}  & 
    \frac{4 \delta k ^2 }{4 \delta ^2 k ^2 + 1} & 
   - \frac{2 \delta k ^3 \beta }{4 \delta ^2 k ^2 + 1} \hat{u} ^{\mathrm{R}} 
   + \frac{ k^2 \beta }{4 \delta ^2 k ^2 + 1} \hat{u} ^{\mathrm{I}} \\[1.3em]
 0 & 0 & \frac{1}{ l_c }
\end{array} \left.\rule{0cm}{2.3cm}\right] 
\end{equation*}
and
\begin{equation*} \bm{ D } =
\left[\rule{0cm}{1.4cm}\right. \begin{array}{c}
 0 \\
 0 \\
 \frac{1}{ \sqrt{ l_c } }
\end{array} \left.\rule{0cm}{1.4cm}\right] \; .
\end{equation*}
The eigenvalues of $\bm{\gamma}$ are 
\begin{equation*}
 \lambda_{1,2} = \frac{4 \delta k ^2 }{4 \delta ^2 k ^2 + 1} \pm \frac{ 2 k }{4 \delta ^2 k ^2 + 1} \imag
 \hspace{15pt} \mathrm{and} \hspace{15pt} \lambda _3 = \frac{1}{ l_c } \; .
\end{equation*}
Thus, since all of the eigenvalues of $\bm{\gamma}$ have positive real parts (recall that $\delta > 0$), there exists a stationary 
distribution for the SDE system (\ref{SDE L0}) and hence there exists an integrable nontrivial solution to (\ref{stationaryFP}) (see \cite{Risken}).  At this point it becomes apparent why 
it is necessary to introduce the parameter $\delta$ as we did in (\ref{regularization}) in order to take the limit $\epsilon \rightarrow 0$.  The 
parameter $\delta$ introduces (arbitrarily small) damping into (\ref{SDE L0}), which 
results in the existence of a stationary distribution for this SDE system.  If $\delta = 0$ then the eigenvalues $\lambda _1$ and $\lambda _2$ of $\bm{\gamma}$ 
are purely imaginary and equation (\ref{stationaryFP}) does not have an integrable nontrivial solution.

We take $\rho ^{\infty}$ to be the normalized solution of (\ref{stationaryFP}), i.e., the solution such that
\begin{equation*}
 \int _{\mathbb{R} ^3} \rho ^{\infty} \diff \hat{v} ^{\mathrm{R}} \diff \hat{v} ^{\mathrm{I}} \diff \eta = 1 \; .
\end{equation*}
Since $\rho _0$ does not depend on $\hat{v} ^{\mathrm{R}} , \hat{v} ^{\mathrm{I}} ,$ or $\eta$, the solvability condition (\ref{solvability condition}) is then equivalent to
\begin{align*}
 \frac{\partialn \rho _0}{\partialn z} = &
 \int _{\mathbb{R} ^3} \rho ^{\infty} \left(\rule{0cm}{.9cm}\right.  \left[\rule{0cm}{.8cm}\right. \hat{v} ^{\mathrm{R}} 
                 \frac{\partialn }{\partialn \hat{u} ^{\mathrm{R}}} + \hat{v} ^{\mathrm{I}} \frac{\partialn }{\partialn \hat{u} ^{\mathrm{I}}} 
  + \frac{L \ell}{ L_x ^2 } \Bigg( \frac{ 1 }{4 \delta ^2 k ^2 + 1} \Bigg) 
                                  | \bm{\kappa} | ^2 \hat{u}^{\mathrm{R}} \frac{\partialn }{\partialn \hat{v} ^{\mathrm{R}} } \\
 & \hspace{15pt} + \frac{L \ell}{ L_x ^2 } \Bigg( \frac{ 2 \delta k }{4 \delta ^2 k ^2 + 1} \Bigg) 
                                  | \bm{\kappa} | ^2 \hat{u}^{\mathrm{I}} \frac{\partialn }{\partialn \hat{v} ^{\mathrm{R}} }
  + \frac{L \ell}{ L_x ^2 } \Bigg( \frac{ 1 }{4 \delta ^2 k ^2 + 1} \Bigg) 
                                  | \bm{\kappa} | ^2 \hat{u} ^{\mathrm{I}} \frac{\partialn }{\partialn \hat{v} ^{\mathrm{I}} } \\
 & \hspace{15pt} - \frac{L \ell}{ L_x ^2 } \Bigg( \frac{ 2 \delta k }{4 \delta ^2 k ^2 + 1} \Bigg) 
                                  | \bm{\kappa} | ^2 \hat{u} ^{\mathrm{R}} \frac{\partialn }{\partialn \hat{v} ^{\mathrm{I}} }  \left.\rule{0cm}{.8cm}\right]
 \left[\rule{0cm}{.8cm}\right. - \frac{1}{2 k } \hat{v} ^{\mathrm{R}} \frac{\partialn \rho _0}{\partialn \hat{u} ^{\mathrm{I}}} 
             + \delta \hat{v} ^{\mathrm{I}} \frac{\partialn \rho _0}{\partialn \hat{u} ^{\mathrm{I}}} \\
& \hspace{15pt} + \frac{ k  \beta l_c }{2} \hat{u} ^{\mathrm{R}} \eta \frac{\partialn \rho _0}{\partialn \hat{u} ^{\mathrm{I}}} 
  + \frac{1}{2 k} \hat{v} ^{\mathrm{I}} \frac{\partialn \rho _0}{\partialn \hat{u} ^{\mathrm{R}}} + \delta \hat{v} ^{\mathrm{R}} \frac{\partialn \rho _0}{\partialn \hat{u} ^{\mathrm{R}}}
  - \frac{ k \beta l_c }{2} \hat{u} ^{\mathrm{I}} \eta \frac{\partialn \rho _0}{\partialn \hat{u} ^{\mathrm{R}}} \left.\rule{0cm}{.8cm}\right] 
    \left.\rule{0cm}{.9cm}\right) \diff \hat{v} ^{\mathrm{R}} \diff \hat{v} ^{\mathrm{I}} \diff \eta 
\end{align*}
which after some simplifying becomes
\begin{align}
\label{solvability condition 2}
  \frac{\partialn \rho _0}{\partialn z} = &
   - \frac{1}{2 k } \frac{\partialn ^2 \rho _0}{ \partialn \hat{u} ^{\mathrm{R}} \partialn \hat{u} ^{\mathrm{I}}} 
                 \int _{\mathbb{R} ^3} \left( \hat{v} ^{\mathrm{R}} \right) ^2 \rho ^{\infty} \diff \hat{v} ^{\mathrm{R}} \diff \hat{v} ^{\mathrm{I}} \diff \eta 
 + 2 \delta \frac{\partialn ^2 \rho _0}{ \partialn \hat{u} ^{\mathrm{R}} \partialn \hat{u} ^{\mathrm{I}}} 
               \int _{\mathbb{R} ^3} \hat{v} ^{\mathrm{R}} \hat{v} ^{\mathrm{I}} \rho ^{\infty} \diff \hat{v} ^{\mathrm{R}} \diff \hat{v} ^{\mathrm{I}} \diff \eta \notag \\
& + \frac{ k  \beta l_c }{2} \Bigg[ \frac{\partialn }{\partialn \hat{u} ^{\mathrm{R}}} 
             \left( \hat{u} ^{\mathrm{R}} \frac{\partialn \rho _0}{\partialn \hat{u} ^{\mathrm{I}}} \right) \Bigg]
               \int _{\mathbb{R} ^3} \hat{v} ^{\mathrm{R}} \eta \rho ^{\infty} \diff \hat{v} ^{\mathrm{R}} \diff \hat{v} ^{\mathrm{I}} \diff \eta \notag \\
 & + \frac{1}{2 k} \frac{\partialn ^2 \rho _0}{\partialn \left( \hat{u} ^{\mathrm{R}} \right) ^2} 
                      \int _{\mathbb{R} ^3} \hat{v} ^{\mathrm{R}} \hat{v} ^{\mathrm{I}} \rho ^{\infty} \diff \hat{v} ^{\mathrm{R}} \diff \hat{v} ^{\mathrm{I}} \diff \eta 
  + \delta \frac{\partialn ^2 \rho _0}{\partialn \left( \hat{u} ^{\mathrm{R}} \right) ^2} 
                \int _{\mathbb{R} ^3} \left( \hat{v} ^{\mathrm{R}} \right) ^2 \rho ^{\infty} \diff \hat{v} ^{\mathrm{R}} \diff \hat{v} ^{\mathrm{I}} \diff \eta \notag \\
 &  - \frac{  k \beta l_c }{2} \hat{u} ^{\mathrm{I}} \frac{\partialn ^2 \rho _0}{\partialn \left( \hat{u} ^{\mathrm{R}} \right) ^2}
               \int _{\mathbb{R} ^3} \hat{v} ^{\mathrm{R}} \eta \rho ^{\infty} \diff \hat{v} ^{\mathrm{R}} \diff \hat{v} ^{\mathrm{I}} \diff \eta 
    - \frac{1}{2 k } \frac{\partialn ^2 \rho _0}{\partialn \left( \hat{u} ^{\mathrm{I}} \right) ^2 } 
                             \int _{\mathbb{R} ^3} \hat{v} ^{\mathrm{I}} \hat{v} ^{\mathrm{R}} \rho ^{\infty} \diff \hat{v} ^{\mathrm{R}} \diff \hat{v} ^{\mathrm{I}} \diff \eta \notag \\
 & + \delta \frac{\partialn ^2 \rho _0}{\partialn \left( \hat{u} ^{\mathrm{I}} \right) ^2 } 
                \int _{\mathbb{R} ^3} \left( \hat{v} ^{\mathrm{I}} \right) ^2 \rho ^{\infty} \diff \hat{v} ^{\mathrm{R}} \diff \hat{v} ^{\mathrm{I}} \diff \eta 
  + \frac{ k \beta l_c }{2} \hat{u} ^{\mathrm{R}} \frac{\partialn ^2 \rho _0}{\partialn \left( \hat{u} ^{\mathrm{I}} \right) ^2}
                        \int _{\mathbb{R} ^3} \hat{v} ^{\mathrm{I}} \eta \rho ^{\infty} \diff \hat{v} ^{\mathrm{R}} \diff \hat{v} ^{\mathrm{I}} \diff \eta \notag \\
 &  + \frac{1}{2 k} \frac{\partialn ^2 \rho _0}{\partialn \hat{u} ^{\mathrm{I}} \partialn \hat{u} ^{\mathrm{R}}} 
                            \int _{\mathbb{R} ^3} \left( \hat{v} ^{\mathrm{I}} \right) ^2 \rho ^{\infty} \diff \hat{v} ^{\mathrm{R}} \diff \hat{v} ^{\mathrm{I}} \diff \eta \notag \\
  &  - \frac{ k \beta l_c }{2} \Bigg[ \frac{\partialn }{\partialn \hat{u} ^{\mathrm{I}}} 
       \left( \hat{u} ^{\mathrm{I}} \frac{\partialn \rho _0}{\partialn \hat{u} ^{\mathrm{R}}} \right) \Bigg]
         \int _{\mathbb{R} ^3} \hat{v} ^{\mathrm{I}} \eta \rho ^{\infty} \diff \hat{v} ^{\mathrm{R}} \diff \hat{v} ^{\mathrm{I}} \diff \eta \notag \\
  &  - \frac{L \ell }{2 L_x ^2 k } | \bm{\kappa} | ^2 \hat{u}^{\mathrm{R}} \frac{\partialn \rho _0}{\partialn \hat{u} ^{\mathrm{I}}}
             + \frac{L \ell }{2 L_x ^2 k} | \bm{\kappa} | ^2 \hat{u} ^{\mathrm{I}} \frac{\partialn \rho _0}{\partialn \hat{u} ^{\mathrm{R}}}  \; .
\end{align}

The normalized solution $\rho ^{\infty}$ of (\ref{stationaryFP}) is the stationary density for (\ref{SDE L0}), which is a 
mean-zero Gaussian \cite{Risken}.  Thus, the values of the integrals in equation (\ref{solvability condition 2}) are 
equal to the corresponding entries of the covariance matrix of the Gaussian $\rho ^{\infty}$.  The components of this covariance matrix satisfy a 
linear system of equations whose coefficients are determined by the 
coefficients in (\ref{stationaryFP}) (see \cite{Risken}).  By solving this system we find the values of the integrals in
(\ref{solvability condition 2}), which are listed in 
Appendix A. 
The fact that equation (\ref{stationaryFP}) does not have an integrable nontrivial solution when $\delta = 0$ can be seen in the 
observation that the components of the covariance matrix of $\rho ^{\infty}$ that are given by (\ref{vR variance}) and (\ref{vI variance})
go to infinity as $\delta \rightarrow 0$.  

We substitute the values of the integrals given in Appendix A
into (\ref{solvability condition 2}) and simplify the resulting equation to get
\begin{align*}
  \frac{\partialn \rho _0}{\partialn z} = 
   & \; \frac{L \ell }{2 L_x ^2  k} | \bm{\kappa} | ^2 \hat{u} ^{\mathrm{I}} \frac{\partialn \rho _0}{\partialn \hat{u} ^{\mathrm{R}}} 
     - \frac{L \ell }{2 L_x ^2  k } | \bm{\kappa} | ^2 \hat{u}^{\mathrm{R}} \frac{\partialn \rho _0}{\partialn \hat{u} ^{\mathrm{I}}} 
     - \frac{ k^3 \beta ^2 l_c 
             \Big( 2 k \hat{u} ^{\mathrm{R}} \big( 1 + \delta l_c ^{-1} \big) - l_c ^{-1} \hat{u}^{\mathrm{I}} \Big) }
              {4 \left( l_c ^{-2} + 4 k^2 \big( 1 + \delta l_c ^{-1} \big) ^2 \right)} \frac{\partialn \rho _0}{\partialn \hat{u} ^{\mathrm{R}}} \\
 & - \frac{ k^3 \beta ^2 l_c
             \Big( 2 k \hat{u} ^{\mathrm{I}} \big( 1 + \delta l_c ^{-1} \big) + l_c ^{-1} \hat{u}^{\mathrm{R}} \Big) }
              {4 \left( l_c ^{-2} + 4 k^2 \big( 1 + \delta l_c ^{-1} \big) ^2 \right)} \frac{\partialn \rho _0}{\partialn \hat{u} ^{\mathrm{I}}} 
              + \frac{ k^2 \beta ^2 l_c }{8} 
              \left( \hat{u}^{\mathrm{I}} \right) ^2 \frac{\partialn ^2 \rho _0}{\partialn \left( \hat{u} ^{\mathrm{R}} \right) ^2} \\
 &    - \frac{ k ^2 \beta ^2 l_c }{4} \hat{u}^{\mathrm{R}} \hat{u}^{\mathrm{I}} 
              \frac{\partialn ^2 \rho _0}{ \partialn \hat{u} ^{\mathrm{R}} \partialn \hat{u} ^{\mathrm{I}}} 
     + \frac{ k^2 \beta ^2 l_c }{8} 
             \left( \hat{u}^{\mathrm{R}} \right) ^2 \frac{\partialn ^2 \rho _0}{\partialn \left( \hat{u} ^{\mathrm{I}} \right) ^2 } \; .
\end{align*}
This is the limiting backward Kolmogorov equation as $\epsilon \rightarrow 0$.  The corresponding (It\^o) SDE system is
\begin{align*}
 \diff \hat{u} ^{\mathrm{R}} _z = & \; \frac{L \ell }{2 L_x ^2 k  } | \bm{ \kappa } | ^2 \hat{u} ^{\mathrm{I}} _z \diff z 
        - \frac{ k^3 \beta ^2 l_c
             \Big( 2 k \hat{u} ^{\mathrm{R}} _z \big( 1 + \delta l_c ^{-1} \big) - l_c ^{-1} \hat{u}^{\mathrm{I}} _z \Big) }
                {4 \left( l_c ^{-2} + 4 k^2 \big( 1 + \delta l_c ^{-1} \big) ^2 \right)} \diff z
                - \frac{ k \beta \sqrt{l_c} }{2} \hat{u}^{\mathrm{I}} _z \diff W _z \\
 \diff \hat{u} ^{\mathrm{I}} _z = & - \frac{L \ell }{2 L_x ^2 k } | \bm{ \kappa } | ^2 \hat{u} ^{\mathrm{R}} _z \diff z 
       - \frac{ k^3 \beta ^2 l_c  
            \Big( 2 k \hat{u} ^{\mathrm{I}} _z \big( 1 + \delta l_c ^{-1} \big) + l_c ^{-1} \hat{u}^{\mathrm{R}} _z \Big) }
             {4 \left( l_c ^{-2} + 4 k^2 \big( 1 + \delta l_c ^{-1} \big) ^2 \right)} \diff z 
               + \frac{ k \beta \sqrt{l_c} }{2} \hat{u}^{\mathrm{R}} _z \diff W _z \; .
\end{align*}
By expressing the above system in terms of the complex state variable $\hat{u} = \hat{u} ^{\mathrm{R}} + \imag \hat{u} ^{\mathrm{I}}$, and using the 
definition (\ref{Fresnel number}) of the Fresnel number $N_{ \mathrm{F}}$, we get
\begin{align*}
 \diff \hat{u} _z = & - \frac{ \imag }{4 \pi N_{ \mathrm{F}}} | \bm{ \kappa } | ^2 \hat{u} _z \diff z 
 - \frac{ k^3 \beta ^2 l_c
             \Big( 2 k \big( 1 + \delta l_c ^{-1} \big) + \imag l_c ^{-1} \Big) }
                {4 \left( l_c ^{-2} + 4 k^2 \big( 1 + \delta l_c ^{-1} \big)^2 \right)} \hat{u} _z \diff z 
                + \imag \frac{ k \beta \sqrt{ l_c }}{2} \hat{u} _z \diff W _z 
\end{align*}
or, equivalently,
\begin{align*}
 \diff \hat{u} _z = & - \frac{ \imag }{4 \pi N_{ \mathrm{F}}} | \bm{ \kappa } | ^2 \hat{u} _z \diff z 
 - \frac{ k^2 \beta ^2 l_c }
                {8 \bigg(    1 + l_c ^{-1} \Big( \delta - \frac{ \imag }{2 k} \Big) \bigg) } \hat{u} _z \diff z 
                  + \imag \frac{k \beta \sqrt{ l_c } }{2} \hat{u} _z \diff W _z \; .
\end{align*}
Taking the inverse Fourier transform (with respect to (\ref{Fourier transform})) to return to the lateral spatial coordinates then gives equation (\ref{first limiting equation}):
\begin{align*}
 \diff u (z, \bm{x}) = & \; \frac{ \imag }{4 \pi N_{ \mathrm{F}}} \mathrm{\Delta} _{\bm{x}} u (z, \bm{x}) \diff z  \\
 & - \frac{ k^2 \beta ^2 l_c }
                {8 \bigg( 1 + l_c ^{-1} \Big( \delta - \frac{ \imag }{2 k} \Big) \bigg) } u (z, \bm{x}) \diff z + \imag \frac{k \beta \sqrt{l_c} }{2} u (z, \bm{x}) \diff W (z) \; .
\end{align*}
Finally, by taking the regularization parameter $\delta$ to zero we obtain the simultaneous paraxial white-noise approximation (\ref{limiting equation}):
\begin{align}
\label{limiting equation rewritten}
 \diff u (z, \bm{x}) = & \; \frac{ \imag }{4 \pi N_{ \mathrm{F}}} \mathrm{\Delta} _{\bm{x}} u (z, \bm{x}) \diff z  \notag \\
 & - \frac{ k^2 \beta ^2 l_c }
                {8 \Big( 1 - \imag \frac{1}{2 k l_c } \Big) } u (z, \bm{x}) \diff z  + \imag \frac{k \beta \sqrt{l_c} }{2} u (z, \bm{x}) \diff W (z) \; .  \hspace{5pt}
\end{align}


\section{Discussion}

We have studied the simultaneous paraxial white-noise limit of the Helmholtz equation in the regime in which the wavelength is of the same order as the correlation length of the refractive index fluctuations in the main direction of propagation.  We showed that this simultaneous limit can be taken in this regime by introducing an arbitrarily small
regularization parameter $\delta$ as we did in equation (\ref{regularization}).  We then took this parameter $\delta$ to zero after taking the simultaneous paraxial white-noise limit in order to obtain the paraxial white-noise approximation (\ref{limiting equation}) for this regime.  We showed that in this regime there is a coupling between the paraxial and white-noise limits that is due to an interaction between the wavelength and the correlation length of the medium fluctuations.  

We now discuss the relationship between the 
paraxial white-noise approximation (\ref{limiting equation}) for the regime in which the wavelength is of the same order as the correlation length of the medium fluctuations and the widely used paraxial white-noise approximation for the high-frequency regime where the wavelength is small compared to the medium correlation length \cite{Bailly,Garnier2009,Garnier2009WM}.  We define the dimensionless parameter $\mu = \lambda _0 / (2 \pi \ell _c )$ and express the paraxial white-noise approximation
(\ref{limiting equation}) in terms of $\mu$:
\begin{align}
\label{limiting equation Discussion mu}
  {\rm d} u (z, \bm{x}) =  & \; \frac{{\rm i}}{4 \pi N_{ \mathrm{F}}} \mathrm{\Delta} _{\bm{x}} u (z, \bm{x}) {\rm d}z \notag \\
 & - \Bigg( \frac{ k^2  \beta ^2 l_c }{8} \Bigg)
    \Bigg( \frac{1}{ 1 - {\rm i} \frac{ \mu }{2} } \Bigg)
 u (z, \bm{x}) {\rm d}z + {\rm i} \frac{ k \beta \sqrt{l_c} }{2}  u (z, \bm{x}) {\rm d}W (z)  \; .  \hspace{2pt}
\end{align}
We expand the coefficient of the second term on the right-hand side in terms of the dimensionless parameter $\mu$:
\begin{equation}
\label{noise induced drift mu}
   - \Bigg( \frac{ k^2  \beta ^2 l_c }{8} \Bigg)
    \Bigg( \frac{1}{ 1 - {\rm i} \frac{ \mu }{2} } \Bigg) =
          - \frac{ k^2 \beta ^2 l_c }{8}
            - {\rm i} \frac{ \mu k^2 \beta ^2 l_c }{16} + ... \; .
\end{equation}
Substituting the first term on 
the right-hand side of (\ref{noise induced drift mu}) in (\ref{limiting equation Discussion mu}) gives the previously derived paraxial
white-noise approximation for the high-frequency regime 
$\mu = \lambda _0 / (2 \pi \ell _c ) \ll 1$.  The second term on the right-hand
side of (\ref{noise induced drift mu}) gives the next-order correction in this 
high-frequency regime.  However, this second 
term as well as subsequent terms in the expansion (\ref{noise induced drift mu}) are order-one in the regime $\lambda _0 \sim \ell _c$.  In this way,
equation (\ref{limiting equation}) can be viewed as a generalization of the previously derived paraxial white-noise approximation for the high-frequency regime that extends the range of validity to include the regime in which the wavelength is of the same order as the correlation length of the medium fluctuations.

One immediate insight that can be gained from the paraxial white-noise approximation (\ref{limiting equation}) is the decay rate of the mean, or coherent, field.  The coherent field decays exponentially with propagation distance in the $z$ direction as scattering by the random medium fluctuations results in energy being transferred to the incoherent field.
The equation for the coherent field $U = E[u]$ is obtained by taking 
the expectation of (\ref{limiting equation}):
\begin{align}
\label{equation for coherent field}
 \diff U (z, \bm{x}) = & \; \frac{ \imag }{4 \pi N_{ \mathrm{F}}} \mathrm{\Delta} _{\bm{x}} U (z, \bm{x}) \diff z 
 - \frac{ k^2 \beta ^2 l_c }
                {8 \Big( 1 - \imag \frac{1}{2 k l_c } \Big) } U (z, \bm{x}) \diff z \; . \hspace{10pt}
\end{align}
It follows from (\ref{equation for coherent field}) that the exponential decay constant $\Lambda$ for the coherent field is given by 
\begin{equation}
\label{decay constant}
 \Lambda = \mathrm{Re} \left(\rule{0cm}{.8cm}\right.  \frac{ k^2 \beta ^2 l_c }
                {8 \Big(    1 - \imag \frac{1}{2 k l_c } \Big) }  \left.\rule{0cm}{.8cm}\right) 
                    = \frac{ k^2 \beta ^2 l_c }{ 8 \left( 1 + \frac{1}{4 k^2 l_c ^2 } \right) } \; .
\end{equation}
In the high-frequency regime $\lambda _0 \ll \ell _c$, the decay constant (\ref{decay constant}) is approximately equal to, but 
slightly smaller than, the one given by the previously derived paraxial white-noise approximation \cite{Garnier2015moment}.
 
Since the correlation length of the refractive index fluctuations due to turbulence varies substantially over different regions of the atmosphere, the propagation regime we have considered in this article describes numerous different atmospheric propagation scenarios. 
As we have mentioned, this correlation length can range from meters near the ground to a few kilometers in the free atmosphere, and can be 
as large as 10 kilometers in the ionosphere.
Thus, the new paraxial white-noise approximation we have derived, and its implications such as the decay rate (\ref{decay constant}) of the coherent field, are relevant for microwave and radiowave propagation through various regions of the atmosphere.  

%

\section*{Appendix A: Covariance matrix of stationary distribution}
\label{Appendix A}
Here we describe the stationary density $\rho ^{\infty}$ corresponding to the SDE system (\ref{SDE L0}), i.e., the normalized solution of the stationary 
Fokker-Planck equation (\ref{stationaryFP}).  The stationary density for (\ref{SDE L0}) is a mean-zero Gaussian \cite{Risken} and therefore it is completely 
determined by its covariance matrix.  The components of this covariance matrix satisfy a 
linear system of equations whose coefficients are determined by the 
coefficients in (\ref{stationaryFP}) (see \cite{Risken}).  By solving this system we find the covariance matrix of $\rho ^{\infty}$, which 
we list below componentwise.  The entries of the covariance matrix of $\rho ^{\infty}$ are given by
\begin{subequations}
\begin{align}
\label{vR variance}
  \int _{\mathbb{R} ^3} & \left( \hat{v} ^{\mathrm{R}} \right) ^2 \rho ^{\infty} \diff \hat{v} ^{\mathrm{R}} \diff \hat{v} ^{\mathrm{I}} \diff \eta = \notag \hspace{230pt} \\
& \hspace{10pt} \frac{ k^2 \beta ^2 \bigg( l_c ^{-1} \left( \hat{u}^{\mathrm{R}} \right) ^2
 + 4 \delta l_c ^{-1} k \hat{u}^{\mathrm{R}} \hat{u} ^{\mathrm{I}} + \Big( l_c ^{-1} + 8 \delta k^2 \big( 1 + \delta l_c ^{-1} \big) \Big) \left( \hat{u}^{\mathrm{I}} \right) ^2 \bigg)}
 {16 \delta \left( l_c ^{-2} + 4 k^2 \big( 1 + \delta l_c ^{-1} \big)^2 \right)} \hspace{2pt}
\end{align}
\begin{align}
  \int _{\mathbb{R} ^3} \hat{v} ^{\mathrm{R}} \hat{v} ^{\mathrm{I}} \rho ^{\infty} \diff \hat{v} ^{\mathrm{R}} \diff \hat{v} ^{\mathrm{I}} \diff \eta = & \;
 \frac{ k^3 \beta ^2 \left( l_c ^{-1} \left( \hat{u}^{\mathrm{I}} \right) ^2
 - 4 k \hat{u}^{\mathrm{R}} \hat{u} ^{\mathrm{I}} \big( 1 + \delta l_c ^{-1} \big) - l_c ^{-1} \left( \hat{u}^{\mathrm{R}} \right) ^2 \right)}
 {8 \left( l_c ^{-2} + 4 k^2 \big( 1 + \delta l_c ^{-1} \big)^2 \right)} \hspace{2pt}
\end{align}
\begin{align}
\label{vI variance}
  \int _{\mathbb{R} ^3} & \left( \hat{v} ^{\mathrm{I}} \right) ^2 \rho ^{\infty} \diff \hat{v} ^{\mathrm{R}} \diff \hat{v} ^{\mathrm{I}} \diff \eta = \notag \hspace{230pt} \\
& \hspace{10pt} \frac{ k^2 \beta ^2 \bigg( l_c ^{-1} \left( \hat{u}^{\mathrm{I}} \right) ^2
 - 4 \delta l_c ^{-1} k \hat{u}^{\mathrm{R}} \hat{u} ^{\mathrm{I}} + \Big( l_c ^{-1} + 8 \delta k^2 \big( 1 + \delta l_c ^{-1} \big) \Big) \left( \hat{u}^{\mathrm{R}} \right) ^2 \bigg) }
 {16 \delta \left( l_c ^{-2} + 4 k^2 \big( 1 + \delta l_c ^{-1} \big)^2 \right)} \hspace{2pt} 
\end{align}
\begin{align}
& \int _{\mathbb{R} ^3} \hat{v} ^{\mathrm{R}} \eta \rho ^{\infty} \diff \hat{v} ^{\mathrm{R}} \diff \hat{v} ^{\mathrm{I}} \diff \eta = 
 - \frac{ k^2 \beta \Big( 2 k \hat{u} ^{\mathrm{I}} \big( 1 + \delta l_c ^{-1} \big) + l_c ^{-1} \hat{u}^{\mathrm{R}} \Big) }
 {2 \left( l_c ^{-2} + 4 k^2 \big( 1 + \delta l_c ^{-1} \big)^2 \right)} \hspace{5pt}
\end{align}
\begin{align}
& \int _{\mathbb{R} ^3} \hat{v} ^{\mathrm{I}} \eta \rho ^{\infty} \diff \hat{v} ^{\mathrm{R}} \diff \hat{v} ^{\mathrm{I}} \diff \eta = 
\frac{ k^2 \beta \Big( 2 k \hat{u} ^{\mathrm{R}} \big( 1 + \delta l_c ^{-1} \big) - l_c ^{-1} \hat{u}^{\mathrm{I}} \Big) }
 {2 \left( l_c ^{-2} + 4 k^2 \big( 1 + \delta l_c ^{-1} \big)^2 \right)} \hspace{5pt}
\end{align}
and
\begin{equation}
 \int _{\mathbb{R} ^3} \eta ^2 \rho ^{\infty} \diff \hat{v} ^{\mathrm{R}} \diff \hat{v} ^{\mathrm{I}} \diff \eta = \frac{1}{2} \; . 
    \hspace{5pt}
\end{equation}
\end{subequations}


\end{document}